\documentclass[11pt,english]{article}

\usepackage{graphicx}
\usepackage{amsmath}
\usepackage{amssymb}
\usepackage{wasysym}
\usepackage{latexsym}
\usepackage{crop}
\usepackage{algorithm,algcompatible}
\usepackage{graphicx}
\usepackage{caption,subcaption}
\usepackage{multirow}
\usepackage{multicol} 
\usepackage{bm}
\usepackage{bbm}
\usepackage{enumerate}
\usepackage{framed} 
\usepackage[framed]{ntheorem}
\usepackage{url}
\usepackage{array}
\usepackage{paralist}
\usepackage{authblk}
\usepackage{booktabs}
\usepackage{sidecap}
\usepackage{wrapfig}
\usepackage{times}
\usepackage{subfiles}
\usepackage{tikz}
\usetikzlibrary{positioning}

\usepackage[sort,numbers]{natbib}
\usepackage{enumitem}
\setlist[enumerate]{leftmargin=*,wide=0em, itemsep=2pt, parsep = 1.5pt, topsep=0pt, label = {\bfseries \arabic*.},listparindent=1.5em}
\setlist[itemize]{leftmargin=*,wide=0em,itemsep=-1pt,topsep=0pt}

\usepackage{titlesec}
\titleformat*{\section}{\normalsize\bfseries}
\titleformat*{\subsection}{\normalsize\bfseries}
\titleformat*{\subsubsection}{\normalsize\bfseries}
\titleformat*{\paragraph}{\normalsize\bfseries}
\titleformat*{\subparagraph}{\normalsize\bfseries}

\usepackage[colorlinks = true, pdfstartview = FitV, linkcolor = blue, citecolor = blue, urlcolor = blue]{hyperref}

\usepackage[capitalise]{cleveref}
\crefname{equation}{}{}
\crefname{figure}{Figure}{Figure}
\creflabelformat{equation}{\textup{(#2#1#3)}}

\usepackage{xspace}

\usepackage{algorithm,algcompatible}

\usepackage{arydshln}
\usepackage{pifont}
\newtheorem{problem}{Problem}

\newcommand{\myqed}{\hskip 10pt $\blacksquare$}

\begin{document}
	\title{Large Markov Decision Processes \\ and Combinatorial Optimization}
\author{
	Ali Eshragh\thanks{International Computer Science Institute, Berkeley, CA, USA.}
}	
\date{}
\maketitle

\begin{abstract}
	Markov decision processes continue to gain in popularity for modeling a wide range of applications ranging from analysis of supply chains and queuing networks to cognitive science and control of autonomous vehicles. Nonetheless, they tend to become numerically intractable as the size of the model grows fast. Recent works use machine learning techniques to overcome this crucial issue, but with no convergence guarantee. This note provides a brief overview of literature on solving large Markov decision processes, and exploiting them to solve important combinatorial optimization problems.
\end{abstract}

\section{Introduction}

Stochastic processes are mathematical models for analyzing the behavior of dynamic systems in the presence of uncertainty. They have a wide range of applications spanning from supply chains and energy systems to epidemiology and computational complexity with applications in big data problems. One common theoretical and practical applications of stochastic processes is in \emph{Markov decision processes}.

Markov decision processes (MDPs) have proved their effectiveness and advantages in modeling and optimizing stochastic dynamic systems evolving in time \cite{Puterman_MDP_2005}. They continue to gain in popularity for modeling a wide range of applications spanning from analyzing and optimizing supply chains and queuing networks to cognitive science and controlling autonomous vehicles (e.g., see \cite{Saghafian_MSOM_2014,Boloori_MSOM_2020,Sennott_SDPCQS_1998,Levine_IJRR_2018,White_I_1985,Carr_ACC_2018,Powell_EJoTL_2012,Brock_SARSOP_2009,Apostolopoulos_CS&PH_2020}). During the past few years, there has been a rapid growth in utilizing MDPs in conjunction with the state-of-the-art machine learning techniques and deep neural networks in solving very large complex stochastic optimization problems. Two excellent recent examples are AlphaGo \cite{Silver_N_2016} for the game of Go, achieving a $99.8\%$ winning rate against other Go programs, and defeating the human European Go champion by $5$ games to $0$, and AlphaZero \cite{Silver_arXiv_2017} for the game of chess, playing much better than all chess programs and convincingly defeating a world-champion. \par

It is evident that MDPs can be applied to a very broad range of optimization problems from deterministic to stochastic, from combinatorial optimization to optimal control, and from one decision maker to multiple--player games \cite{Filar_CMDP_1996}. Nonetheless, they are plagued by a crucial practical issue, the so-called \emph{curse of dimensionality}, implying that the size of these models grows exponentially fast and consequently, finding an ``exact'' optimal solution for such large models becomes intractable in terms of both time and storage complexity. For instance, the size of the state space for an MDP corresponding to the game of Tetris is more than the total number of all molecules in the world \cite{Gabillon_NIPS_2013}! 

To overcome the curse of dimensionality, two main approaches have been developed: (i) approximation, and (ii) simulation. While the former utilizes different architectures to reduce the dimension of the large MDP (e.g., via aggregation \cite{Bertsekas_IEEE_2019}) or develops an approximate version of the exact algorithms (e.g., approximate dynamic programming \cite{Powell_ADP_2007}), the latter uses a computer model in place of a mathematical model (e.g., simulation-based algorithms \cite{Chang_SAMDP_2013}) or recent advances in deep neural networks (e.g., deep reinforcement learning \cite{Sutton_RL_2018}).  

\section{State Aggregation Method} 

This is a fast growing approach to approximate large MDPs which is the main focus of this note
\vspace*{-0.15cm}
\setlist[itemize,1]{leftmargin=*,wide=0em, noitemsep,nolistsep}
\begin{algorithm}[htbp]
	\caption{Generic Framework of State Aggregation Method}
	\begin{algorithmic}[1]
		\STATEx \textbf{Input:} A large Markov decision process 
		\begin{description}
			\item[\textit{Step 1.}] Find similar states and aggregate them together. 
			\item[\textit{Step 2.}] Construct an aggregated MDP based on the aggregated states. 
			\item[\textit{Step 3.}] Solve the aggregated MDP exactly and find an optimal policy. 
			\item[\textit{Step 4.}] Disaggregate the optimal policy by mapping to the original MDP. 
		\end{description}
		\vspace{1mm}
		\STATEx \textbf{Output:} The disaggregated policy
	\end{algorithmic}
	\label{alg:typical_state_agg}
\end{algorithm}
\vspace*{-0.15cm}
\hspace*{-0.4cm}The generic framework of the method, depicted in \cref{alg:typical_state_agg}, seems sensible except that the proverbial ``devil is in the details''. How, precisely, should Steps $1-4$ be implemented to achieve a provably fast and robust algorithm? Most recent works use machine learning techniques to establish features for distinguishing similarities between states \cite{Bertsekas_IEEE_2019}. More specifically, a neural network with these features as input is trained to provide a set of aggregated states as output. Although, the numerical experiments have revealed promising results, they are not guaranteed to work for all or even most problems, as there is no theoretical error bounds or convergence results on these approximation methods \cite{Bertsekas_RLOC_2020}. Of course, there exists a few other state aggregation methods with convergence guarantee for some specific classes of MDPs (e.g., see \cite{Hutter_TCS_2016})), but they mainly use the results from theoretical computer science and do not exploit the rich literature of MDPs. To address this crucial limitation of the current methods, we recently developed a new state aggregation framework based on the theory of MDPs (rather than a machine learning method), and tested this new framework on several numerical examples. Preliminary results are very promising.  

\section{Combinatorial Optimization} 

Combinatorial optimization (CO) is an important branch of mathematics with a broad range of applications in operations research (e.g., see \cite{Applegate_TSP_2007,Mazyavkina_arXiv_2020,Li_JCO_2013,Arlotto_OR_2014}). Arguably, one of the central problems in CO is the traveling salesman problem (TSP) as several other CO problems are either a variation of TSP, or can be transformed to TSP. Although TSP is an NP-hard problem, that is, no one has found a solution algorithm with polynomial worst case runtime for it, yet, its simple appearance and wide range of applications have captivated many researchers. Much of its difficulty is concealed in its combinatorial structure, which is characterized by its NP-complete counterpart problem, namely the Hamiltonian cycle problem (HCP).

In $1994$, CI Filar, together with Krass, \cite{Filar_MOR_1994} developed a new model for HCP by embedding this problem in a perturbed MDP. One significant benefit of this stochastic approach is that it creates an opportunity to search for Hamiltonian cycles among extreme points of a suitably constructed, polyhedral subset of the polytope induced by constrained MDPs.
Continuing this line of research, Feinberg \cite{Feinberg_MOR_2000} and later Eshragh et al. have improved these polyhedral domains and exposed their geometrical and structural properties. Moreover, Eshragh and Filar postulated some important conjectures concerning the complexity of certain random walk algorithms on the extreme points of these polytopes induced by constrained MDPs. These conjectures which have been further refined very recently, will provide a new approach to tackle not only HCP, but also its optimization variant, TSP.

\section{Nexus between MDP and CO} 

One may suggest to develop theoretically sound frameworks for the state aggregation method for solving large MDPs, and exploit them to solve TSP by embedding it into suitably constrained MDP. To achieve these two aims, the following two problems should be addressed: 
\begin{problem}\label{obj:MDP}
	Develop a new theoretical state aggregation framework for solving large MDPs and provide theoretical error bound as well as guaranteed convergence.
\end{problem}
\begin{problem}\label{obj:TSP}
	Develop a new solution approach for TSP by embedding it into a suitably constrained MDP and utilizing results of \cref{obj:MDP}.
\end{problem}

\bibliographystyle{abbrvnat}
\bibliography{./bib_Ali}

\end{document}